%% file: Highly_Scalable_arxiv.tex
\documentclass[hidelinks,onefignum,onetabnum]{siamart250106}

\newcommand{\mycomment}[1]{%

\ifthenelse{\isodd{\value{page}}}{%
\makeatletter\normalmarginpar%
\marginpar{\tiny {#1}}%
}{%
\makeatletter\reversemarginpar%
\marginpar{\tiny {#1}}%
}}%

\input{ex_shared}

\ifpdf
\hypersetup{
  pdftitle={Highly Scalable Two-level Monolithic Overlapping Schwarz Preconditioners for Thermo-elastoplastic laser beam welding problems},
  pdfauthor={T. Bevilacqua, A. Klawonn and M. Lanser}
}
\fi

\begin{document}

\maketitle
\begin{abstract}
A thermo-elastoplastic finite element approach is used to perform the simulation of a Laser Beam Welding (LBW) process. This results in a nonlinear, nonsymmetric saddle point multiphysics system, for which the nonlinearity is handled via the Newton method. The iterative solution of the arising linear system is accelerated by using robust and highly scalable, overlapping  Schwarz domain decomposition preconditioners. It is well-known that a one-level method of this type is not scalable and therefore a second level has to be added.

Therefore, the construction and numerical analysis of  monolithic, two-level overlapping Schwarz preconditioners with different variants of the GDSW (Generalized Dryja-Smith-Widlund) coarse space are presented here. A new and parallel efficient implementation of several variants of GDSW, that is, GDSW, RGDSW, and GDSW*, in PETSc, is introduced, which is usable for multiphysics problems, as, for example, the thermo-mechanical LBW problem considered here. Different combinations of the GDSW variants for the different fields (temperature and displacements) are compared and parallel scalability for realistic problems is shown.
\end{abstract}



\section{Introduction}
Overlapping Schwarz domain decomposition methods are a class of preconditioners for the large linear systems obtained from the discretization of many linear and nonlinear partial differential equations. Originally, they were introduced for uniformly elliptic partial differential equations. In \cite{klawonn1998overlapping} monolithic overlapping Schwarz methods were introduced for saddle point problems originating from the discretization of Stokes' equations discretized with mixed finite elements. Here, the overlapping domain decomposition is applied directly to the saddle point problem without reducing it to sub-blocks; this requires the solution of smaller saddle point problems on each subdomain. Additionally, a coarse problem based on a saddle point formulation was introduced in \cite{klawonn1998overlapping}. For a brief review of monolithic overlapping Schwarz preconditioners, see the discussion further below in this introduction. In the present paper, we develop a monolithic overlapping Schwarz preconditioner for a challenging multiphysics problem from \ac{LBW}.

\ac{LBW} is a noncontact joining process that recently has become considerably more important in the course of the increasing degree of automation in industrial production. The advantages of this process are its short cycle times and the small heat-affected zones. However, the high cooling rate of the metal specimen inherent in the process can lead to solidification cracks, due to a residual melt oversaturated with certain alloy elements. It is therefore crucial to try to quantitatively understand the development mechanisms of these fractures and their correlation to process parameters. In particular, this is the focus of the DFG research group 5134 ''Solidification Cracks During Laser Beam Welding -- High Performance Computing for High Performance Processing''\footnote[3]{https://www.for5134.science/en/}, where experiments, advanced multiscale finite element modelling, and scientific computing are brought together to obtain new insights. A nice example for the interplay of experiments and high-performance simulations within this project is presented in \cite{tp145laser}, where a \ac{CTW} test is performed both experimentally and virtually using finite element simulations. The \ac{CTW} test is used to check if a material is susceptible for cracks during welding by measuring its ductility in the mushy zone, i.e. the region where a material is partially solid and partially liquid during cooling and phase transition. The probability of cracks is increased in the \ac{CTW} test by applying different strain conditions and then temperature and strain are tracked optically to see how they affect cracking. This is a crucial analysis to improve welding processes and prevent defects. In~\cite{tp145laser}, our PETSc~\cite{balay2019petsc}-based software library FE2TI~\cite{klawonn2020computational,nakajima} and its numerical solver infrastructure was used to perform the finite element simulations of the \ac{CTW} test on high-performance computers. In the present work here we provide an accompanying in-depth study of the performance and scalability of our FE2TI software for \ac{LBW} problems and, to be more precise, we investigate the robustness and scalability of the \ac{DDMs} used to solve the large finite element problems arising in \ac{LBW} simulations.

The thermomechanical properties of \ac{LBW} can be represented by a coupled system of time dependent thermo-elastoplasticity equations \cite{SimMie:1992:act}. Due to the nature of the problem, the linear system that arises from a finite element discretization in space, a Backward Euler discretization in time, and a linearization by Newton's method is ill-conditioned. As a consequence, here we discuss robust and parallel scalable preconditioners to accelerate the convergence of the \ac{GMRES} method. Since we are predominantly interested in analyzing the performance of the domain decomposition preconditioners and not in a precise modelling of realistic laser simulations, we here consider finite element problems where the action of the laser beam is simplified. In particular, we assume that the melting pool, that is, the portion of the material molted by the laser, has the shape of a simple cylinder that goes vertically through the entire material. Actually, the exact form of the melting pool does not play an important role for the solver performance and simulations of thermo-elastic and thermo-elastoplastic \ac{LBW} processes with more accurate representation of the melting pool \cite{bakir2018numerical} can be found in \cite{tp145laser} and \cite{bevilacqua2024monolithic}.

The action of the laser is applied as a volumetric constraint to the temperature field in the region of the melting pool that moves over time with a problem dependent speed. More realistic heat source models for the laser can also used in order to obtain more accurate simulations \cite{hartwig2024physically,hartwig2023numerical,hartwig2023volumetric}, but  again we do not expect a significant influence on the performance of our preconditioners.

As already mentioned, in this work, we focus on the \ac{DDMs} used to iteratively solve the linearized systems arising in the simulation of \ac{LBW} processes. To be more precise, we make use of a monolithic overlapping two-level additive Schwarz preconditioner for saddle point problems. Such monolithic  domain decomposition preconditioners have been first introduced in \cite{klawonn1998overlapping}; see also \cite{klawonn2000comparison} for a comparison to block-diagonal and block-triangular preconditioners using overlapping Schwarz methods on the sub-blocks.
In order to obtain numerical scalability, a coarse level needs to be introduced. 
 In \cite{klawonn1998overlapping} Lagrangian coarse spaces that operate on the saddle point problem defined on the discretization introduced by the domain decomposition have been used. Such a construction of the coarse space using the domain decomposition as a coarse finite element triangulation is quite restrictive. To obtain more flexibility in the construction of the coarse space, especially for more general shapes of subdomains, we use variants of the Generalized Dryja-Smith-Widlund (GDSW) coarse spaces, in which coarse basis functions are defined using problem dependent extension operators. In the case of symmetric positive definite matrices these extensions are discrete harmonic  or, more generally, energy minimizing. The GDSW coarse space was originally introduced in \cite{dohrmann2008domain,dohrmann2008family} for scalar diffusion problems and has been  extended to saddle point problems arising from the discretization of incompressible fluid flow problems in \cite{heinlein2019monolithic,hochmuth2020parallel}. Let us note that in \cite{heinlein2019monolithic} a fully monolithic approach was chosen, using saddle point problems to compute the extensions.
Despite being very robust, a disadvantage of this type of coarse space is the rapid growth of its size when scaling to a large number of subdomains, especially for three dimensional problem. Several strategies to reduce the size of the coarse space have been introduced, for example, GDSW*~\cite{hochmuth2020parallel} and Reduced-GDSW~\cite{dohrmann2017design} coarse spaces. In the case of multi-physics problems one can combine the different coarse spaces to reduce its overall size. For instance, in case of fluids, choosing a larger coarse space for the velocity block and a smaller one for the pressure space is often benefitial. This idea has been introduced and analyzed in ~\cite{sassmannshausen2025} for Navier-Stokes problems where many combinations of coarse spaces have been considered. Parallel implementations of monolithic Schwarz preconditioners have also been applied to fully-coupled nonlinear chemo-mechanics problems; see~\cite{rheinbach:kiefer} for details.

In our present work, for the first time, we give an exhaustive computational analysis of similar combinations of different coarse spaces of GDSW-type for thermo-elastoplasticity problems. For this purpose, we implemented all variants, that are, GDSW, GDSW*, and RGDSW, in  PETSc. For the nonlinear and time-dependent coupled thermo-elastoplasticity problems multiple linear systems are arising from Newton's method. To deliver a fair comparison and sound results, we also integrated well-known strategies to improve the efficiency, such as, an adaptive time stepping approach, truncation of the coarse matrix, and the possibility to recycle certain solver components from previous Newton steps. We show efficient parallel scalability of our new solvers for complex and realistic \ac{LBW} problems.

\indent The remainder of the paper is organized as follows. In \Cref{sec:thermo-eq} we present the thermo-elastoplasticity equations, in \Cref{sec:FEM} we introduce the discretization, in \Cref{sec:GDSW} we present the domain decomposition preconditioners and all coarse spaces, in \Cref{sec:impl} we describe the implementation in PETSc, and, finally, in \Cref{sec:results} we show several numerical results and a detailed comparison of the performance of different monolithic coarse spaces.

\section{Thermo-elastoplasticity equations}\label{sec:thermo-eq}
In our work we focus on studying a coupled system of thermo-elasto-plasticity equations
\begin{align*}
 {\rm div} ({\Bsigma(\bm{u},\theta)}) &= 0, \\
 {\rm div} (\bm{q})+\gamma \, {\rm div} (\dot{\bm{u}}) \, \theta + \rho \, c_\rho \, \dot{\theta} &= 0.
\end{align*}
Here $\bm{u}(x,t)$ and $\theta(x,t)$ represent the displacement vector and the absolute temperature respectively. 
These equations are simplified versions of the ones that can be found in \cite{SimMie:1992:act}, since in our case we assume to not have any volume acceleration and no external heat source. The action of a laser beam  is modeled by a volumetric constraint, directly fixing the temperature to the melting temperature of the material in the melting pool, which is an unusual but still effective approach.\\
The other variables, that represent the material-specific and temperature-dependent material parameters are: $c_\rho$ the heat capacity, $\bm{q}$ the heat flux vector, and $\gamma = 3\alpha_T\kappa$ the stress temperature modulus. Moreover, $\alpha_T$ denotes the coefficient of linear thermal expansion, $\kappa = E / 3(1-2\nu)$ the bulk modulus and $E$ the Young's modulus. These temperature-dependent parameters have been chosen conformely to model an austenitic chrome-nickel steel(1.4301) as shown in \Cref{tab:mat_param}. Finally, $\rho = 7.919\, kg/m^3$ is the density of the material, assumed to be constant~\cite{richter2010physical}.\\
Here, the thermo-elastoplastic model assumes an additive decomposition of total strains $\Bve = [\ve_{11} \,\, \ve_{22} \,\, \ve_{33} \,\, \ve_{13} \,\, \ve_{12} \,\, \ve_{23}]$ (in Voigt notation), into elastic, plastic, and thermal components. Based on this, the free energy function $\psi$ is formulated in terms of the thermo-elastic strains, temperature, and an internal variable, leading to the constitutive equation $\Bsigma = \partial_{\Bve^{te}} \psi$, which defines the stress-strain relationship, with $\Bve^{te}$ the thermo-elastic component. The model also considers the balance of momentum and energy and to characterize the plastic behavior, a von Mises yield criterion is employed, incorporating the deviatoric stress dev$\bm{\Bsigma}$. For further details related to the model, see \cite{SimMie:1992:act,SimHug:1998:cin} and for the same specific problem \cite{tp145laser}.

The weak formulation of the coupled thermo-elastoplasticity equations is
\begin{equation}\label{eq:weak}
\begin{split}
 \int_\Omega \Bsigma (\bm{u},\theta)\, :\, \Bve (\bm{v}) \diff x &= 0 \qquad \forall \bm{v} \in [H^1(\Omega)]^3, \\
\int_\Omega \bm{q} \cdot \nabla q \diff x + \int_\Omega \{- 3\alpha_T \, \kappa \,{\rm tr}[\dot{\Bve}(\bm{u})] \, \theta - \rho \,\textit{c}_\rho \dot{\theta} \}q \diff x &= 0 \qquad \,\forall q \in H^1(\Omega),
\end{split}
\end{equation}
with Dirichlet boundary conditions $\bm{u}:=(u_x,u_y,u_z) = \bm{u}_D$ on $\partial \Omega_{D,u} \subset \partial \Omega$. As already mentioned, the only heat entry is caused by the laser beam, which we implement as a volume of molten mateial with fixed melting temperature $\theta_l$. Therefore we have the additional constraint $\theta := \theta_{l} \gg~\!\!0$ in $\Omega_{\theta} \subset \overline{\Omega}$. Let us note that the displacement is not fixed within $\Omega_{\theta}$ and the material can still move. This is a nonlinear and nonsymmetric saddle point system for which the theoretical framework is still an object of studies, for example, an investigation of pairs of inf-sup-stable finite elements. A theoretical investigation of this problem is beyond the scope of this work.

\begin{table*}[t]
     \caption{Parameters of the material austenitic chrome-nickel steel(1.4301) at different temperatures. The values, taken from~\cite{hartwig2023volumetric}, have been provided by Bundesanstalt f{\" u}r Materialforschung und -pr{\" u}fung (BAM).  A linear interpolation is used between the temperature intervals.}\label{tab:mat_param}
  \begin{tabular}{|l@{\hskip 0.0cm}|@{\hskip 0.1cm}c@{\hskip 0.05cm}|@{\hskip 0.1cm}c@{\hskip 0.05cm}|@{\hskip 0.1cm}c@{\hskip 0.05cm}|@{\hskip 0.1cm}c@{\hskip 0.05cm}|@{\hskip 0.1cm}c@{\hskip 0.05cm}|@{\hskip 0.1cm}c@{\hskip 0.05cm}|@{\hskip 0.1cm}c@{\hskip 0.05cm}|@{\hskip 0.1cm}c@{\hskip 0.05cm}|@{\hskip 0.1cm}c@{\hskip 0.05cm}|}
    \hline
    \multirow{2}{*}{Parameter} & \multicolumn{9}{c|}{Value} \\
                         & \multicolumn{9}{c|}{(Temperature $[C]$)}\\
    \hline
    \multirow{2}{*}{$E \cdot 1e4 \, [{\rm N/mm}^2]$} & 20.0 & 19.1 & 17.5 & 12.5 & 7.2 & 1.6 & 0.1 & -- & --\\
     & (20) & (170) & (400) & (800) & (1000) & (1100) & (1500) & -- & -- \\[1mm]     
    \multirow{2}{*}{$\nu $} & 0.271 & 0.284 & 0.300 & 0.319 & 0.329 & 0.364 & 0.364 & -- & -- \\ 
     & (20) & (183) & (484) & (799) & (994) & (1994) & (2000) & -- & -- \\[1mm]
    \multirow{2}{*}{$\alpha_T \cdot 1e-5 \, [{\rm K}^{-1}]$} & 1.6 & 1.81 & 1.98 & 2.13 & 2.23 & 2.23 & 2.33 & -- & -- \\ 
    & (20) & (200) & (580) & (1000) & (1200) & (1500) & (2000) & -- & -- \\[1mm]
    \multirow{2}{*}{$\bm{\lambda} [{\rm W}/({\rm m}\cdot {\rm K}))]$} & 15.6 & 18.1  & 21.0 & 23.8 & 26.6  & 34.4 & 35.0 & 60 & -- \\
    & (20) & (200) & (400) & (600) & (800) & (1350) & (1393) & (1460) & -- \\[1mm]    
    \multirow{2}{*}{$c_{\rho}  \cdot 1e5 \, [{\rm J/kg}\cdot {\rm K}]$ }& 5.11 & 5.42 & 5.75 & 6.05 & 6.30 & 6.85 & 7.30 & 20.20 & 50.00 \\
    & (20) & (200) & (400) & (600) & (800) & (1350) & (1427) & (1442) & (1460) \\
  \hline
\end{tabular}
\end{table*}

\section{Finite element discretization}\label{sec:FEM} 
The equations in \eqref{eq:weak} are linearized by applying the Newton method to the nonlinear boundary value problem. This yields, in each nonlinear iteration, a linear problem which spatially is discretized by using a mixed finite element method; see \cite{hartwig2024physically} for further details on the linearization. Next, we provide details on the mixed finite element discretization. 

Let $\tau_h$ be a uniform mesh of $N_e$ hexahedral elements $\Omega_e$ of $\Omega$ with characteristic mesh size $h$. We introduce the conforming discrete piecewise linear displacement and temperature spaces
\begin{align*}
V^h &= V^h(\Omega) = \{\bm{u}\in [\mathcal{C}^0(\Bar{\Omega}) \cap H^1(\Omega)]^3: \bm{u}|_T \in Q_1 \ \forall T \in \tau_h\}, \\
Q^h &= Q^h(\Omega) = \{\theta \in \mathcal{C}^0(\Bar{\Omega}) \cap H^1(\Omega): \theta|_T \in Q_1 \ \forall T \in \tau_h\}
\end{align*}
respectively, of Q1-Q1 mixed finite elements, where the displacement space is discretized using the $B_{bar}$ formulation to prevent the phenomenon of volumetric locking \cite{SimTayPis:1985:vap}. Here, $\mathcal{C}^0(\Bar{\Omega})$ denotes the space of continuous functions on $\Bar{\Omega}$ and $H^1(\Omega)$ the usual Sobolev space. Since we do not have a stable theoretical framework, we can not ensure that this choice of finite elements satisfies an inf-sup condition \cite{boffi2013mixed}. As a time discretization we use a Backward Euler method with time step $\Delta t$. Therefore, in each time step, we have to solve a linearized system until the global absolute residual of Newton's method falls below a fixed tolerance. The resulting discrete saddle point problem has the form
\begin{equation*}
    K\varDelta \mathbf{d} = \begin{bmatrix}
        K_{uu} & K_{u\theta} \\
        K_{\theta u} & K_{\theta \theta}
    \end{bmatrix} \begin{bmatrix}
        \varDelta \mathbf{d}_u \\ \varDelta \mathbf{d}_\theta 
    \end{bmatrix} = \begin{bmatrix}
        R_u \\ R_\theta 
    \end{bmatrix} = R, 
\end{equation*}
where $\varDelta \mathbf{d}_u$ and $\varDelta \mathbf{d}_\theta$ represent the Newton update for the displacement and temperature, and $R_u$ and $R_\theta$ the vectors of the residual, respectively. The block matrices are obtained by finite element assembly and we have
\begin{equation*}
\begin{split}
& K_{uu} = \mathbb{A}^{N_e}_{e=1} \bigg[ \int_{\Omega_e} \mathbf{B}_u^T \, \mathbb{C}\, \mathbf{B}_u \diff x \bigg], \qquad K_{u\theta} = \mathbb{A}^{N_e}_{e=1} \bigg[ -\int_{\Omega_e} \mathbf{B}_u^T \, (3\alpha_T \kappa \mathbf{1})\, {\rm \mathbf{N}}_\theta \diff x \bigg],\\
& \hspace{20mm} K_{\theta u} = \mathbb{A}^{N_e}_{e=1} \bigg[-\frac{1}{\varDelta t} \int_{\Omega_e} \theta \, {\rm \mathbf{N}}_\theta^T \, (3\alpha_T \kappa \mathbf{1}^T) \, \mathbf{B}_u \diff x \bigg], \\ 
& K_{\theta\theta} = \mathbb{A}^{N_e}_{e=1} \bigg[ -\int_{\Omega_e} \mathbf{B}_\theta^T \, \bm{\lambda} \, \mathbf{B}_\theta \diff x \, - \int_{\Omega_e} {\rm \mathbf{N}}_\theta^T \, {\rm tr}[\dot{\Bve}(\bar{\bm{u}})] \, {\rm \mathbf{N}}_\theta \diff x \,- \frac{1}{\varDelta t}\int_{\Omega_e} {\rm \mathbf{N}}_\theta^T \,\rho \, c_\rho \, {\rm \mathbf{N}}_\theta \diff x\bigg],
\end{split}
\end{equation*}
where $\mathbf{N}_u$ and ${\rm \mathbf{N}}_\theta$ are the finite element nodal basis functions for displacement and temperature, $\mathbf{B}_u$ and $\mathbf{B}_\theta$ denote their derivatives in tensor notation, and $\bar{\bm{u}}$ are the displacement components resulting from  the previous Newton iteration.
The nature of this linear system requires suitable preconditioners to accelerate the convergence of an iterative solver, as, in our case, the \ac{GMRES} method. Our choice is to use two-level monolithic Schwarz type preconditioners.

\section{Monolithic Overlapping Schwarz Preconditioners}\label{sec:GDSW}

Here, we exclusively consider Schwarz \ac{DDMs} to precondition the linearized systems. They have the advantage that they can be used to directly precondition the original system without any modifications and it is additionally possible to add monolithic second levels or coarse spaces to obtain robustness and scalability.
\subsection{Monolithic two-level Schwarz preconditioner with GDSW coarse space}\label{sec:GDSWc}
More precisely, we consider a monolithic overlapping Schwarz domain decomposition preconditioner of GDSW (Generalized Dryja-Smith-Widlund) type for the saddle point problem. This is a two-level Schwarz preconditioner where the coarse space and the local problems have the same block structure as the original one. We refer the interested reader to \cite{dohrmann2008domain,dohrmann2008family} for a complete definition of the GDSW preconditioner for elliptic problems and to  \cite{heinlein2019monolithic} for saddle point problems resulting from Stokes' equations.

We introduce $\{ \Omega_i\}^N_{i=1}$, a nonoverlapping decomposition of the domain $\Omega$ into $N$ subdomains of characteristic diameter $H$, and $\{ \Omega'_i\}^N_{i=1}$, the corresponding overlapping domain decomposition with $k$ finite element layers of overlap, generated by geometric or graph techniques. We define the interface $\Gamma$ of the nonoverlapping domain decomposition
$\Gamma=\{ x\in (\Omega_i \cap \Omega_j ) \setminus \partial \Omega_D : i \neq j, 1\leq i, j\leq N \},$ that is divided into $M$ connected components, specifically into $M_V$ vertices, $M_E$ edges, and $M_F$ faces. Here, all interface nodes belonging to the same two subdomains build a face and all nodes belonging to the same set of more than just two subdomains build an edge. Finally, edges which consist of just one node are considered as vertices. These edges, faces, and vertices are used later on to define the basis functions of the GDSW coarse space. In general, the definition of an interface allows for the partitioning of the \ac{DOFs} into interface ($\Gamma$) and interior ($I$) ones and, after a reordering, an equivalent partitioning of $K$, that is,  
\begin{equation*}
K:= \begin{bmatrix}
        K_{II} & K_{I\Gamma} \\
        K_{\Gamma I} & K_{\Gamma \Gamma}
    \end{bmatrix}.
\end{equation*}
We further introduce the operators $R_{u,i} : V^h \rightarrow V^h_i$ and $R_{\theta,i} : Q^h \rightarrow Q^h_i$ for $i = 1,...,N$. These are restrictions from the global finite element spaces to the local finite element spaces defined on the overlapping subdomains. The monolithic restriction operators are then denoted by $R_i : V^h \times Q^h \rightarrow V^h_i \times Q^h_i$ for $i=1,...,N$, with
\begin{equation*}
R_i:= \begin{bmatrix}
        R_{u,i} & 0 \\
        0 & R_{\theta,i}
    \end{bmatrix}.
\end{equation*} 
The transposed operators $R^T_i$, $R^T_{u,i}$, and $R^T_{\theta,i}$ are the corresponding prolongation operators.
Now, the additive monolithic two-level GDSW preconditioner can be written as
\begin{equation*}
\hat{B}^{-1}_\text{GDSW} = \Phi K^{-1}_0 \Phi^T + \sum_{i=1}^N R_i^T \, K_i^{-1} \, R_i,
\end{equation*}
where the local stiffness matrices $K_i$ are extracted from $K$ by
\begin{equation*}
K_i = R_i K R_i^T , \qquad {\rm for}\,\, i = 1, . . . , N,
\end{equation*} 
and the coarse operator is defined by
\begin{equation}\label{eq:coarse}
K_0=\Phi^T K \Phi.
\end{equation}
Here, $\Phi$ is the matrix that contains the coarse space basis functions $\phi^i$, $i = 1,...,N_\phi$, column-wise. These functions are constructed as the discrete saddle point extensions of functions $\phi^i_\Gamma$ defined on the interface by solving the linear system
\begin{equation}\label{interfprb}
\begin{bmatrix}
        K_{II} & K_{I \Gamma} \\
        0 & I
    \end{bmatrix} \begin{bmatrix}
        \phi^{i}_I \\ \phi^i_{\Gamma}
    \end{bmatrix} = \begin{bmatrix}
        0 \\ \phi^i_{\Gamma}
    \end{bmatrix},
\end{equation}
where we decompose the interface values of the coarse basis into displacements ($u_0$) and temperature ($\theta_0$) as $\phi^i_{\Gamma} = [ \,\phi^{i\,T}_{\Gamma,u_0} \ \phi^{i\,T}_{\Gamma,\theta_0}\,]^T$.
The entries of $\phi^i_{\Gamma}$ are usually set to 1 or 0 and generally represent the restriction of the null space of the operators $K_{uu}$ and $K_{\theta \theta}$ on the interface components  $\Gamma_j$ with $j=1,...,M$. 

In three dimensions the null space of the $K_{uu}$ block is given by the six rigid body modes (three translations and three rotations) defined on the interface components. For the $K_{\theta \theta}$ operator, these consist of the functions that have constant temperature on $\Omega$, so we restrict the constant function 1 to vertices, edges, and faces.

Specifically, each node of our mesh can be associated to a vector with four components, where the first three entries represent the displacements and the last one the temperature. We introduce the vectors
\begin{equation}\label{eq:tra}
    v_{u,1}=\begin{bmatrix}1 \ 0 \ 0 \ 0\end{bmatrix}^T, \quad
    v_{u,2}=\begin{bmatrix}0 \ 1 \ 0 \ 0\end{bmatrix}^T, \quad
    v_{u,3}=\begin{bmatrix}0 \ 0 \ 1 \ 0\end{bmatrix}^T, \quad
    v_{\theta}=\begin{bmatrix}0 \ 0 \ 0 \ 1\end{bmatrix}^T
\end{equation}
where the first three are the translations into $x$-, $y$-, and $z$- directions and the last one is the restriction of the constant temperature on the node. Moreover we introduce the vectors
\begin{equation}\label{eq:rot}
    v_{u,4}=\begin{bmatrix} y \ -x \ 0 \ 0\end{bmatrix}^T, \quad
    v_{u,5}=\begin{bmatrix} -z \ 0 \ x \ 0\end{bmatrix}^T, \quad
    v_{u,6}=\begin{bmatrix} 0 \ z \ -y \ 0\end{bmatrix}^T, \quad
\end{equation}
that represent the three rotation with center placed in the origin $(0,0,0)$.
For each vertex, edge, and face $\Gamma_j,\;j=1,...,M,$ of the interface $\Gamma$, we then construct four basis functions  $\phi_{\Gamma_j}^k$, $k=1,...,4$, by filling the  \ac{DOFs} associated to the nodes of $\Gamma_j$ with $v_{u,k}$ for  $k=1,2,3$ and $v_\theta$ for $k=4$. If rotations are used, we add $\phi_{\Gamma_j}^k$, $k=5,6,7,$ using $v_{u,k-1}$ analogously. For the nodes in $\Gamma \setminus \Gamma_j$, $\Phi_{\Gamma_j}^k$ is set to zero. Finally, we collect all $\Phi_{\Gamma_j}^k,\; j=1,...,M, \; k=1,...,4$, as the columns of $\Phi_\Gamma$ and thus obtain $N_\phi = 4M$ coarse basis functions and columns in $\Phi$. In case of using rotations, we analogously obtain $N_\phi = 7M$ coarse basis functions.

After having solved \Cref{interfprb} for each basis function of the coarse space, the matrix $\Phi$ has the following block representation
\begin{equation}\label{eq:fullbasis}
\Phi = \begin{bmatrix}
        \Phi_{u,u_0} & \Phi_{u,\theta_0} \\
        \Phi_{\theta, u_0} & \Phi_{\theta,\theta_0}
    \end{bmatrix}.
\end{equation}
With this construction, the presence of nonzero blocks on the off-diagonal can deteriorate the numerical scalability, especially if nonzero Dirichlet boundary conditions for the temperature are present (cf. \cite{heinlein2019monolithic}), as we have with the laser beam. Therefore, these blocks need to be removed, obtaining 
\begin{equation}\label{eq:basis}
\Phi = \begin{bmatrix}
        \Phi_{u,u_0} & 0 \\
        0 & \Phi_{\theta,\theta_0}
    \end{bmatrix}.
\end{equation}
Let us remark that the extensions from the interface to the interior parts of the subdomains, that is, solving \Cref{interfprb}, can be done in parallel in a perfectly scalable way, since $K_{II}$ is a block diagonal matrix. The factorization of $K_{II}$ can thus be done block-wise on each subdomain simultaneously using a sparse direct solver. Additionally, it is not necessary to compute the extension to all subdomains for all basis functions. For example, for the four or seven coarse basis functions belonging to a face, only the extension to the two neighboring subdomains are necessary, the extensions to all other subdomains are zero and must not be computed. Consequently, computing $\Phi$ is a parallel scalable process.

\subsection{GDSW* and RGDSW coarse spaces}\label{sec:RGDSW}
In the GDSW coarse space, the number of coarse basis functions depends on the number of components in the equivalence classes into which the interface between subdomains is divided—namely, vertices, edges, and faces. This number can grow rapidly, especially in 3D domain decomposition. To manage this complexity, alternative techniques such as GDSW* and RGDSW have been introduced.

We briefly outline the ideas behind these alternative coarse spaces and refer to \cite{sassmannshausen2025,hochmuth2020parallel} for further details. The key difference lies in how the components of the equivalence classes are treated. In the GDSW* coarse space, unlike the original approach, only two interface components are used; the first one is obtained by each vertex along with its adjacent edges and the latter by faces. Since two adjacent vertices share an edge, the resulting interface components overlap over this edge. In order to maintain the robustness of these coarse spaces a partition of unity of the basis functions defined across the interface, according to \eqref{eq:tra} and \eqref{eq:rot} has to be preserved, therefore an inverse multiplicity scaling (Option 1 in Section 3~\cite{dohrmann2017design}) is applied.\\

Similarly, in the RGDSW coarse space, only one type of interface component is considered, that is each vertex along with its adjacent edges and faces. Again, multiplicity scaling is necessary to handle the overlap between interface components, which in this case happens for edges and faces. Having partitioned the interface into  edges,  vertices, and faces, the coarse space size of GDSW is thus $C\cdot (M_E+M_V+M_F)$. The size of the GDSW* coarse space is smaller with $C\cdot (M_V+M_F)$ and finally the size of the RGDSW coarse space is the smallest one with only $C \cdot M_V$ coarse basis functions. Here, $C$ is either 4 or 7, depending on the usage of rotations.\\ 
Finally, we point out that for a two-dimensional domain decomposition of a three-dimensional problem, the GDSW and GDSW* coarse spaces coincide, as no vertices are available in the decomposition.

\subsection{Combination of coarse spaces}\label{sec:Comb}
In \Cref{sec:GDSWc} and \Cref{sec:RGDSW}, we explained the monolithic approach to construct GDSW, GDSW*, and RGDSW coarse spaces for the complete saddle point problem. However, in this particular instance and for more general multi-field problems, it is feasible to amalgamate these approaches by independently selecting the most suitable one for each field. In fact, the sole requirement is to provide a partition of unity of the basis function across the interface. Examining the construction of these in equations \eqref{eq:tra} and \eqref{eq:rot} reveals that the problems in equation \eqref{interfprb} are independent of each other. \\
Let us give an explicit example for the construction of a GDSW*(T)-RGDSW coarse space, that means to choose GDSW* that includes only translations for the elastoplastic block and RGDSW for the thermal block. First of all, we note that in the matrix $\Phi$ in \eqref{eq:fullbasis} we have, columnwise represented, $N_{u}$ and $N_{\theta}$ basis function respectively associated to the displacements and temperature fields, with $N_{\phi} = N_u + N_{\theta}$. As we saw before, the number of coarse basis functions $N_\phi$ depends on the number of components into which the interface is divided. The choice GDSW*(T)-RGDSW means that the interface is split into $M_V+M_F$ components for the elastoplastic block and into $M_V$ components for the thermal block. Since in this coarse space, only translations are considered, we construct $3$ interface basis functions for each interface component by using the vectors $v_{u,i}, i=1,..,3$ in \eqref{eq:tra} and we then extend them to the subdomain by solving \eqref{interfprb}, obtaining $N_u = 3\cdot (M_V + M_F)$. In a similar way we define the interface basis function for the temperature by only using the vector $v_\theta$ and again we extend the solution to internal part again solving \eqref{interfprb}, obtaining $N_{\theta} = M_V$.

In this way, for instance in our particular problem, one could opt for a more extensive space for the elastoplastic block eventually deciding to introduce also the rotations and a more limited space for the thermal one. These available choices contribute to the flexibility of the preconditioners, enabling their design and adaptation for different problems, for example by leveraging the most economical coarse space when feasible and reserving the more costly one where it is necessary. As we will show in the numerical results, a more extensive selection for the elastoplastic block is required in our particular problem, while the more limited coarse space can be utilized for the temperature.

\section{Implementation}\label{sec:impl}
All Schwarz preconditioners and GDSW-type preconditioners discussed here are incorporated into the FE2TI software package, a C/C++ library that utilizes a hybrid MPI/OpenMP~\cite{chandra2001parallel} parallelization strategy, leveraging the PETSc library~\cite{balay2019petsc}. To enhance its capabilities, we have integrated an interface with the FEAP library~\cite{taylor2014feap}, allowing access to FEAP’s finite element formulations and material models. FE2TI has already demonstrated an excellent strong and weak scalability behavior, efficiently handling simulations with over one million parallel ranks for realistic dual-phase steel deformation scenarios; for further details, see~\cite{klawonn2020computational,nakajima}. In FE2TI, efficient implementations of the domain decomposition approaches BDDC (Balancing Domain Decomposition by Constraints) and FETI-DP (Finite Element Tearing and Interconnecting - Dual Primal) are already integrated (see~\cite{bddc1,klawonn2020computational,bddc2}) and we were able to reuse some of the efficient building blocks for the new Schwarz preconditioners. For example, faces, edges, and vertices of the domain decomposition are identified in a scalable way using the implementation from our BDDC implementation.\\ 
Within the Schwarz preconditioners, our software defaults to the more efficient restricted option~\cite{frommer2001algebraic} with  an overlap of one layer of finite elements for the first level, although the classical additive alternative is also available. Regarding the second (or coarse) level, several options are supported, that are, GDSW, GDSW*, and RGDSW. While this work focuses on solving a thermo-elastoplasticity problem, that is, a two-by-two block system, we note that our coarse solver implementation is also applicable to general multifield problems. In fact, as described in \Cref{sec:Comb}, we can select different approaches for each field, allowing for flexible combinations in the choice of the coarse space to better adapt to the nature of the problem.  We provide a detailed comparison of different combinations in our numerical results.\\
 The coarse basis functions are constructed by solving local problems and then assembled into an MPI parallel matrix $\Phi$ (see \Cref{eq:basis}) shared among all the ranks. During this construction the truncation tolerance is applied by zeroing out small entries of $\Phi$. The coarse matrix $K_0$ (see \eqref{eq:coarse}) is then obtained by using the PtAP product from PETSc. However, the resulting parallel  matrix is spread over all the MPI ranks, which makes its factorization using a sparse direct solver inefficient. Therfore $K_0$ is then redistributed only among a reduced number of MPI ranks within a seperate sub-communicator.
Since the coarse space is typically expensive to compute, we have incorporated the ability to recycle it after its initial computation, enabling a reuse in subsequent applications of the preconditioner across different linear systems. This feature is particularly beneficial in cases such as ours, where a Newton method is applied at each time step, requiring the solution of multiple spectrally similar linear systems.  Different strategies for the recycling can be applied: one can decide to reuse the matrix $\Phi$ and recompute $K_0$ and its factorization, or one can  recycle all, $\Phi$, $K_0$, and its factorization. Clearly, the cheapest option is to reuse both basis functions and coarse matrix, saving an extra factorization. This leads to a slightly worse approximation of the preconditioner that can lead to an increasing number in the \ac{GMRES} iteration in later Newton iterations. Anyway, we will see in the numerical results that this effect is minimal. Therefore, here we only consider the complete recycling.\\
All local problems are solved exactly using the sparse direct solver package PARDISO from the Math Kernel Library (MKL)~\cite{schenk2001pardiso,intel_onemkl_pardiso}, while the coarse problem is solved using MKL Cluster PARDISO (see \Cref{subsec:num:setup}). Therefore, the coarse matrix is distributed to the ranks of an MPI sub-communicator of a reasonable size.\\
In summary, we employ the highly efficient \ac{GMRES} solver from PETSc to iteratively solve large sparse linear systems. To ensure controlled convergence, we use our highly scalable implementation of Schwarz preconditioners with GDSW-type coarse spaces. Additionally, our interface to FEAP enables seamless integration with the models described in~\Cref{sec:FEM}.

\begin{figure}[]
    \includegraphics[width=0.495\textwidth,trim={0 20 50 0},clip]{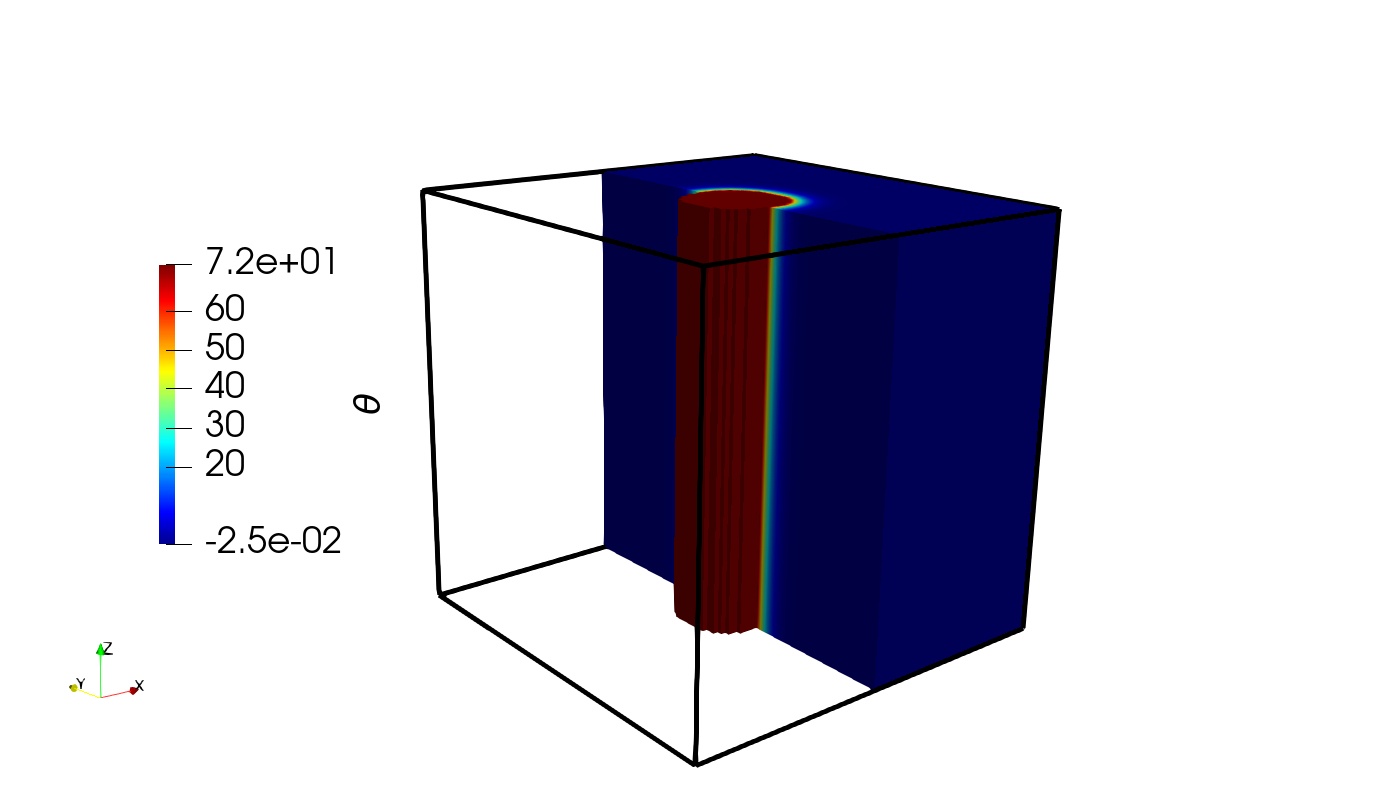}
    \includegraphics[width=0.495\textwidth,trim={0 20 50 0},clip]{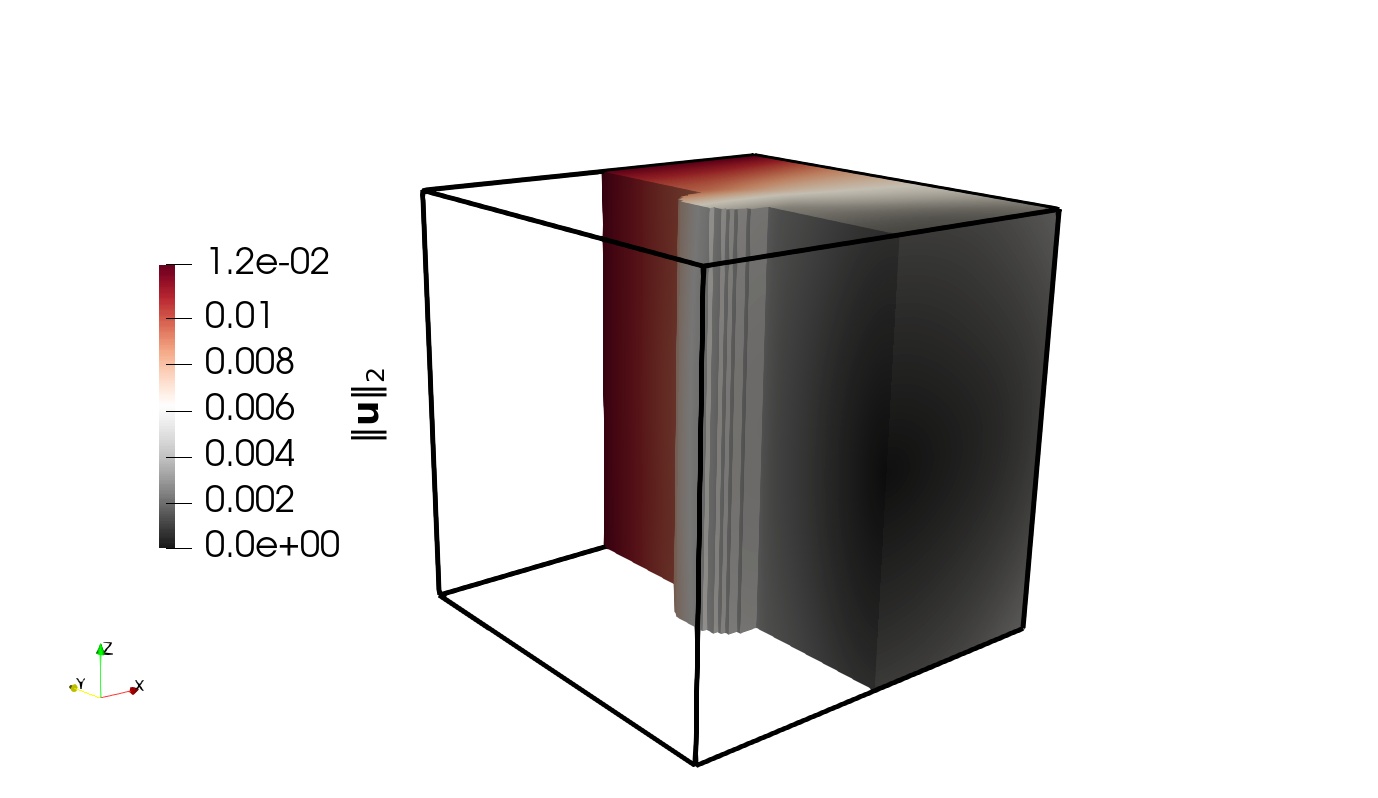}
{
  \caption{Section of temperature field (left) and norm of the displacements (right) for the cube test.}\label{fig:cubetest}
}
\end{figure}

\section{Numerical Results}\label{sec:results}
In our numerical simulations, we consider thermo-elastoplasticity problems that emulate real experiments in which a steel plate of domain $\Omega = [0,l_x] \times[0,l_y] \times[0,l_z]$, with $l_x,l_y,l_z >0$, is melted by a laser beam.  As anticipated in \Cref{sec:thermo-eq}, the effect of the laser is modeled by fixing the temperature field in the melting pool. This is enforced as a gradient with a ratio of $\dot{\theta}_l = 14\,400^\circ C/s$. Specifically, at the beginning of the simulation each DOF is initialized with a temperature $\theta_0 = 20^\circ C$, then in each time step we identify the \ac{DOFs} that are located inside the region heated by the laser, which is modeled as a cylinder of radius $2\,mm$ passing vertically through the entire material.  The temperature of these latter \ac{DOFs} is then increased by a value $\dot{\theta_l} \cdot \Delta t$ until the maximum value of $\theta_l = 1460^\circ C$ is reached. We note that, this heating procedure (initialization phase), where the laser remains fixed in its initial position, takes a total of $0.1$s and will not be completed during the first set of experiments, since only few first time iterations are taken into account. This aspect is relevant for the simulations in \Cref{sec:2Ddom} and \Cref{sec:Appl}, where after the initialization phase the laser moves along the plate.

\begin{figure}[]
    \includegraphics[width=1\textwidth,trim={0 0 0 0},clip]{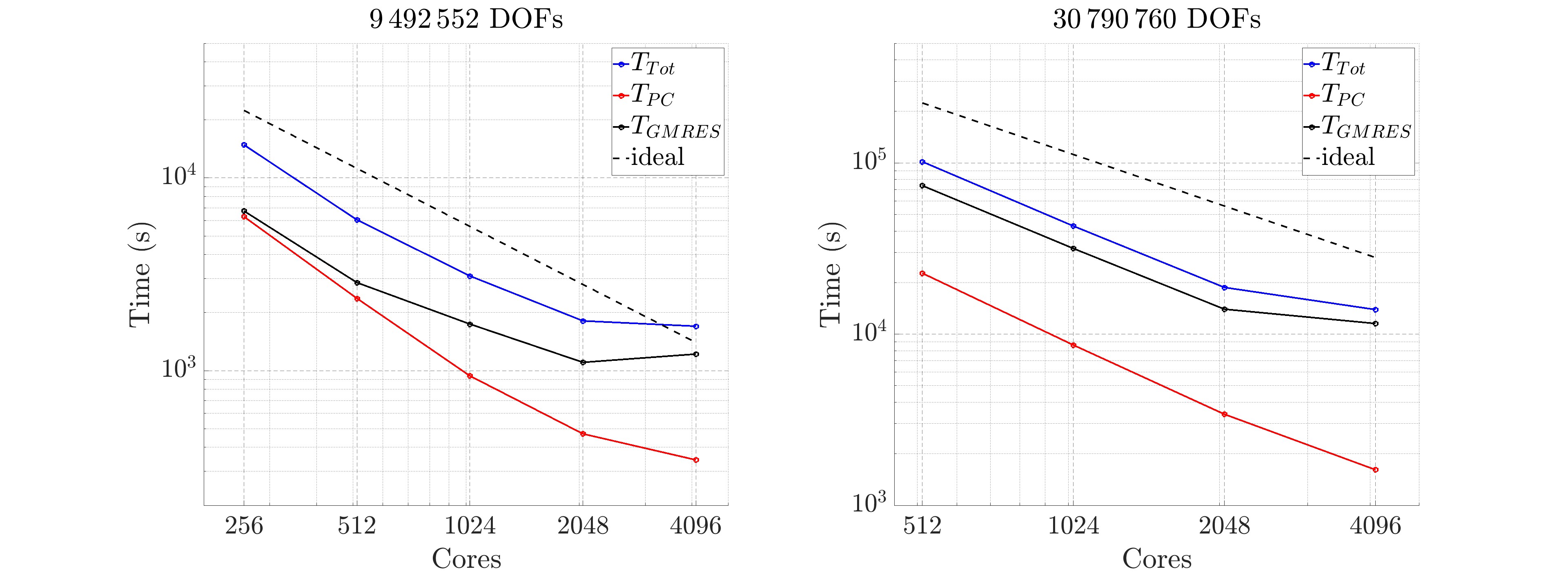}
{
  \caption{Strong scalability for GDSW(T)-RGDSW with recycling option over 130 time iterations.}\label{fig:platetest}
}
\end{figure}

For the displacements, additional Dirichlet boundary conditions are applied to simulate an external load acting on the plate during the welding process. To avoid higher localization of the strain in the corners of the domain, the boundary conditions are applied in the following way: on the surface with $y = 0$, we prevent the movement along the $y$ axis by enforcing $u_y$ = 0 and the rotation of the material by setting $u_x = 0$ along the $z$ axis and $u_z = 0$ along the $x$ axis.  The external load is then applied to the face with $y = l_y$ by enforcing $u_y = u_D(t)$, defined as $u_D(t) := \min(\ve_{22},\dot{\ve}_{22} \cdot t)$ with the external strain $\ve_{22}$ and strain rate $\dot{\ve}_{22}$ specified before the simulations.
We consider two different setups for our simulations. First, to validate our implementation and assess the performance of the solver, we analyze a steel cube subjected to both the laser and an external load using a 3D domain decomposition strategy (\Cref{fig:cubetest}), that is, a decomposition into $N \times N \times N$ cubical subdomains. Once the optimal solver is identified, we apply it to a more realistic case involving a thin steel plate using 2D domain decomposition \Cref{fig:platetest}, that is, a decomposition into $N \times N$ thin subdomains. For an application of our solver to reproduce a physically realistic experiment, we refer the interested reader to \cite{tp145laser}.

 If not explicitely specified, by default all the simulations use a time step size of $\Delta t = 0.001\,\text{s}$ and the extenal strain is applied from the beginning of the simulation with $\ve_{22} = 0.03$ and $\dot{\ve}_{22} = 0.06$. For each time step, Newton's method is applied until the global residual falls below an absolute tolerance of $10^{-4}$, while, for the linear systems, the \ac{GMRES} method is employed, with a stopping criterion of either a relative residual error of $10^{-6}$ or an absolute residual error of $10^{-10}$ based on the unpreconditioned norm.
In our numerical results and tables we report the following statistics: $it_{GMRES}$ is the average count of \ac{GMRES} iteration over all the linear systems, $it_N$ is the average count of Newton iteration per time step, $T_{PC}$ is the total time to assemble the preconditioner,  $T_{Ass}$ is the total time to assemble the linear systems, $T_{Sol}$ is the total time to solve the linear systems, $T_{Tot}$ is the total time for the simulation,  $T_{Navg}$ is the average time per Newton iteration, $T_{Tot}/\#Time it$ is the average time per time iteration, $P_E$ is the parallel efficiency, $Cores$ is the number of cores and MPI ranks used, $Subcomm$ is the number of MPI ranks in the sub-communicator used for the coarse solve, $PC\,Type$ is the type of preconditioner used, $(T+R)$ indicates that in the coarse space for the elastoplasticity field translations and rotations are included while $(T)$ indicates that only translations are considered.

\subsection{General solver setup}\label{subsec:num:setup}
\begin{table}[]
    \centering
    \begin{tabular}{ccccccc}
        Tolerance & $it_{GMRES}$ & $it_{N}$ & $T_{PC}$ & $T_{Sol}$ & $T_{Tot}$  \\
        \hline
        \hline
          0     & 18.9 & 3 & 425.5 & 82.3  & 548.6 \\
          1e-4  & 18.9 & 3 & 326.2 & 58.5  & 416.8  \\
          1e-2  & 42.3 & 3 & 314.3 & 128.8 & 475.4  \\
          \hline
    \end{tabular}
    \caption{ Solver performance with respect to truncation tolerance. $10\times10\times10$ subdomains with $20\times20\times20$ finite elements. All values in $\Phi$ which are lower than the given tolerance are set to zero before building the coarse matrix. All times are given in seconds.}
    \label{tab:tolPhi}
\end{table}
\begin{table}[]
    \centering
    \begin{tabular}{ccccccc}
        $Cores\,(Coarse)$ & \multicolumn{2}{c}{512(10\,465)} & \multicolumn{2}{c}{1\,728(37\,730)} & \multicolumn{2}{c}{4\,096(92\,610)} \\ 
        $Subcomm$ & $T_{PC+Sol}$ & $T_{Tot}$ & $T_{PC+Sol}$ & $T_{Tot}$ & $T_{PC+Sol}$ & $T_{Tot}$ \\
        \hline
        \hline
        1  & 49.5 & 54.4  & 211.2 & 216.9 & x     & x     \\
        4  & 37.7 & 42.4  & 103.1 & 108.3 & 247.3 & 254.9 \\
        8  & 35.5 & 40.2  & 85.1  & 90.3  & 159.5 & 206.0 \\
        16 & 35.7 & 40.4  & 76.7  & 81.8  & 173.7 & 181.8 \\    
        32 & \bf{35.3} & \bf{40.0}  & \bf{68.1}  & \bf{73.5}  & 157.5 & 164.6 \\    
        64 & 37.1 & 42.0  & 85.4  & 96.0  & \bf{150.4} & \bf{157.6} \\    
        \hline
    \end{tabular}
    \caption{ Solver performance with respect to the number of subcommunicator used for the coarse solver. Each subdomain contains $10\times10\times10$ finite elements. All values in $\Phi$ which are lower than $1e-4$ are truncated to zero before building the coarse matrix. All times are given in seconds.}
    \label{tab:subcomm}
\end{table}

 In this section we exlusively consider the simplified \ac{LBW} problem on a cube geometry to find out the optimal setup for our Schwarz preconditioner.
Before comparing and analyzing the performance of different combinations of GDSW-type coarse spaces to identify the optimal one, we first investigate two key aspects of the solver that will remain fixed afterwards. 

One of the main challenges of GDSW-type preconditioners is the rapid growth of the coarse space, which results in a large coarse matrix which tends to be much denser than a typical finite element stiffness matrix. This growth significantly increases the computational cost of the coarse factorization using a sparse direct solver, especially since a forward-backward substitution has to be applied in each \ac{GMRES} iteration. To mitigate this issue, we can adjust the construction of the coarse matrix. Recall that the coarse matrix is given by $K_0=\Phi^T K \Phi$.

A first strategy to improve efficiency is to remove small entries in the coarse basis matrix $\Phi$ that fall below a certain truncation tolerance. This results in a sparser and more computationally efficient coarse matrix. To evaluate this approach, we consider the GDSW*(T+R)-RGDSW coarse space and perform five time iteration steps. In \Cref{tab:tolPhi}, we analyze the solver’s behavior for different threshold values applied to the entries of $\Phi$. Our results indicate that we can be fairly aggressive, using a relatively high tolerance of $1e-4$ without increasing the iteration count. This significantly reduces both the cost of preconditioner construction and the overall solving time. However, as expected, setting the truncation tolerance too high degrades the accuracy of the coarse space, making it a poor approximation.

Additionally, in \Cref{tab:subcomm}, we report the timings for the  same experiment as before fixing a truncation tolerance of $1e-4$ while varying the number of computational cores and the number of MPI ranks in the sub-communicator for the coarse solve. As expected, solving the coarse problem on a single MPI rank is not feasible, especially as the number of MPI processes for the global problem increases, leading to a larger coarse space. Based on the results in \Cref{tab:subcomm}, we determine the optimal number of MPI ranks within the coarse sub-communicator for all further tests.  The best setup for the different number of subdomains is marked in bold face in \Cref{tab:subcomm}.

\subsection{Comparison of different coarse spaces}
\begin{table}[]
    \centering
    \begin{tabular}{ccccccc}
            \hline
            \multicolumn{7}{c}{No recycling}\\
        $PC\,Type$ & $Coarse$ & $it_{GMRES}$ & $it_{N}$ & $T_{PC}$ & $T_{Sol}$ & $T_{Tot}$  \\
        \hline
        \hline
          GDSW(T+R)-GDSW       & 36\,225 &  18.8 & 3 & 481.3 & 68.4  & 582.2 \\
          GDSW(T+R)-RGDSW      & 31\,266  & 18.8 & 3 & 443.7 & 65.1  & 540.6  \\
          GDSW*(T+R)-GDSW*     & 23\,954  & 18.9 & 3 & 338.5 & 60.1  & 430.5  \\
          GDSW*(T+R)-RGDSW     & 21\,303  & 18.8 & 3 & 326.2 & 58.5  & 416.5  \\
          RGDSW(T+R)-RGDSW     & 5\,103   & 48.6 & 3 & 254.2 & 133.1 & 418.9  \\ 
         \hline
         \hline
          GDSW(T)-GDSW       & 23\,265 &  29.4 & 3 & 342.8 & 90.0  & 465.1 \\
          GDSW(T)-RGDSW      & 18\,306  & 29.1 & 3 & 309.2 & 86.8  & 427.6  \\
          GDSW*(T)-GDSW*     & 13\,667  & 29.5 & 3 & 259.6 & 83.4  & 374.7  \\
          GDSW*(T)-RGDSW     & 11\,016  & 29.7 & 3 & 249.2 & 82.5  & 363.2  \\
          RGDSW(T)-RGDSW     & 2\,916   & 53.8 & 3 & 216.7 & 142.9 & 391.4  \\ 
         \hline
    \end{tabular}
    \caption{ Solver performance with respect to different choices for the coarse space. $10\times10\times10$ subdomains with $20\times20\times20$ finite elements. The truncation tolerance is set to $1e-4$. All times are given in seconds.}
    \label{tab:comparison}

    \centering
    \begin{tabular}{ccccccc}
            \hline
            \multicolumn{7}{c}{With recycling}\\
        $PC\,Type$ & $Coarse$ & $it_{GMRES}$ & $it_{N}$ & $T_{PC}$ & $T_{Sol}$ & $T_{Tot}$  \\
        \hline
        \hline
          GDSW(T+R)-GDSW       & 36\,225 &  18.9 & 3 & 230.6 & 68.9  & 331.3 \\
          GDSW(T+R)-RGDSW      & 31\,266  & 18.9 & 3 & 217.8 & 65.5  & 314.8  \\
          GDSW*(T+R)-GDSW*     & 23\,954  & 19.1 & 3 & 183.0 & 60.8  & 275.2  \\
          GDSW*(T+R)-RGDSW     & 21\,303  & 18.9 & 3 & 178.6 & 58.9  & 268.9  \\
          RGDSW(T+R)-RGDSW     & 5\,103   & 48.7 & 3 & 154.4 & 134.3 & 320.1  \\ 
         \hline
         \hline
          GDSW(T)-GDSW       & 23\,265 &  29.8 & 3 & 183.9 & 92.8  & 308.6 \\
          GDSW(T)-RGDSW      & 18\,306  & 29.1 & 3 & 172.9 & 87.3  & 291.5  \\
          GDSW*(T)-GDSW*     & 13\,667  & 29.9 & 3 & 156.2 & 85.0  & 272.5  \\
          GDSW*(T)-RGDSW     & 11\,016  & 29.7 & 3 & 152.6 & 83.3  & 267.0  \\
          RGDSW(T)-RGDSW     & 2\,916   & 53.8 & 3 & 141.9 & 143.8 & 317.2  \\ 
         \hline
    \end{tabular}
    \caption{ Solver performance with respect to different choices for the coarse space with the recycling option active. $10\times10\times10$ subdomains with $20\times20\times20$ finite elements. The truncation tolerance is set to $1e-4$. All times are given in seconds.}
    \label{tab:my_label}
\end{table}

Fixing the optimal basic setup from the previous section, that is, truncating $\Phi$ with a tolerance of $1e-4$ and choosing the optimal size of the coarse sub-communicator, we now compare different coarse spaces.
We again consider five time iteration steps of the thermo-elastoplastic problem on the cube by using $10\times 10\times 10$ cores and subdomains with $20\times 20\times 20$ finite elements each. In \Cref{tab:comparison} we report the effect of different coarse space choices on the solver performance in terms of the average iteration counts in each time step and the total timings.
We recall that the first coarse space is used for the elastoplastic part of the coupled system, while the second one for the temperature, for example, GDSW-RGDSW means that we use GDSW for the mechanics and RGDSW for the temperature. We also investigate the choice to include only translations (T) or both translations and rotations (T+R) in the coarse space for the elastoplastic block.\\
At first sight, we note that, as expected, the time to construct the preconditioner increases with the size of the coarse space. Additionaly, we observe that including rotational modes significantly improves convergence by reducing the number of \ac{GMRES} iterations. This pattern is observed across all configurations, confirming that rotations enhance the preconditioners effectiveness. However, these improvements in iteration counts come at the cost of an increment of the coarse space size, which affects the setup time of the preconditioner as well as the solution time. Among the configurations that include rotations, GDSW*(T+R)-RGDSW provides the best trade-off, maintaining a moderate coarse space size (21\,303) while achieving one of the lowest total solution times (416.5s). On the other hand, RGDSW(T+R)-RGDSW, despite having the smallest coarse space (5\,103), suffers from a significantly higher iteration count (48.6), leading to an overall longer solution time (418.9s), regardless its lower setup cost.

A similar trend is observed in the translation-only (T) configurations. While reducing the coarse space leads to lower preconditioner setup times, the number of iterations increases, impacting the overall efficiency. For example, GDSW*(T)-RGDSW, with a coarse space size of 11\,016, takes 29.7 iterations and has a total computation time of 363.2s, which is an improvement over larger coarse spaces and its (T+R) counterpart, despite having a higher number of iterations. So far, these results indicate that using only translations appear to be a valid option. Especially when reducing the coarse space size has priority. Often, even if we have an increased number of iterations, the total time to solution is lower. Overall, GDSW(T)-RGDSW emerges as the most efficient configuration, achieving a strong balance between coarse space size, iteration count, and total computational time.

\subsection{Recycling option}

In \Cref{tab:my_label}, we further analyze solver performance repeating the previous test but now activating the recycling option for the coarse operator and the coarse basis functions. A clear improvement is observed across all configurations, particularly regarding the preconditioner setup time, which is significantly reduced compared to \Cref{tab:comparison}. We recall that with this approach, in each time step, the coarse space is constructed only once and in the subsequent Newton iterations within the same time step only the sparse direct solvers for the local subdomains have to be updated.
 This strategy mostly affects the setup time for the preconditioner, for instance, with GDSW*-RGDSW (T+R) we have a reduction in the setup time  from 326.2s to 178.6s and at the same time its total solution time drops from 416.5s to 268.9s.  We stress out that, by recycling the coarse space and the basis functions, we are not using anymore the  exact two-level Schwarz preconditioner, but an approximation. This is highlighted by the increase in iteration count (18.8 to 18.9 iterations), that anyway is minimal. This trend is consistent across other configurations, confirming that recycling effectively mitigates the computational overhead associated with the coarse space factorization while maintaining solver robustness. 

 Since, as already mentioned, the most expensive part of the setup of the preconditioner is the construction of the coarse space, the benifits of the recycling are more evident especially when the plastic effect is stronger and, as a result, the number of Newton iterations per time step is increased.  Moreover, since the coarse space is assembled just once, i.e. in the first Newton iteration, the time spent to assemble the preconditioner in the further ones is affected only by the direct solver applied on the subdomain, levelling the timings also for the most expensive coarse space. In that case the time for setting up the preconditioner looses importance, and thus solver setups with fast times the solution get more and more favorable. The results confirm that GDSW*-RGDSW (T+R) remains the optimal configuration, particularly when combined with recycling, as it provides the best trade-off between coarse space size, computational cost, and solver efficiency.

\begin{table}[]
    \centering
    \begin{tabular}{rrrrrrrrr}
        \hline
    	\multicolumn{9}{c}{No recycling}\\
        $Cores$ & $Coarse$ & $it_{GMRES}$ & $it_N$ & $T_{PC}$ & $T_{Sol}$ & $T_{Ass}$ & $T_{Tot}$ & $P_{E}$\\
        \hline
        \hline
        64   & 1\,053  & 29.2 & 7.6  & 10\,898  & 1\,512  & 251.1  & 12\,681 & 1.00    \\
        216  & 4\,115  & 27.4 & 7.6  & 1\,810   & 348.6 & 116.1  & 2\,282& 2.77 \\ 
        512  & 1\,0465 & 27.0 & 7.6  & 517.2  & 139.0 & 50.4   & 708.4 & 2.23 \\
        1728 & 37\,730 & 27.0 & 7.8  & 232.2  & 107.4 & 15.7   & 356.5 & 1.32 \\ 
        4096 & 92\,610 & 26.6 & 7.8  & 213.7  & 126.2 & 7.8    & 349.5 & 0.56 \\
        \hline
    \end{tabular}

    \centering
    \begin{tabular}{ccccccccc}
        \hline    	
    	\multicolumn{9}{c}{With recycling}\\
        $Cores$ & $Coarse$ & $it_{GMRES}$ & $it_N$ & $T_{PC}$ & $T_{Sol}$ & $T_{Ass}$ & $T_{Tot}$ & $P_{E}$ \\
        \hline
        \hline
        64   & 1\,053  & 29.6 & 7.6  & 5\,643  & 1\,519  & 250.1 & 7\,431 & 1.00    \\
        216  & 4\,115  & 27.9 & 7.6  & 810.1 & 356.6 & 115.2 & 1\,288& 2.88 \\ 
        512  & 10\,465 & 27.9 & 7.6  & 213.9 & 144.7 & 49.8  & 409.5 & 2.27 \\
        1728 & 37\,730 & 28.3 & 7.8  & 59.5  & 112.2 & 15.6  & 187.9 & 1.46 \\
        4096 & 92\,610 & 28.6 & 7.8  & 39.1  & 135.2 & 7.8   & 182.2 & 0.63 \\
        \hline
    \end{tabular}
    \caption{Strong scalability for GSDW*(T+R)-RGDSW without recycling option (top) and with recycling option (bottom). $12\,194\,500$ \ac{DOFs}. Computed the first 5 time iterations. All times are given in seconds.}
    \label{tab:strongRe}
\end{table}
\subsection{Strong Scalability}\label{sec:strongCube}
\Cref{tab:strongRe} provide the strong scalability results for the GDSW*-RGDSW coarse space, with and without the recycling option. We performed the test for a mesh of $12\,194\,500$ \ac{DOFs} over the first $5$ time iteration steps. To increase the plasticity effect we apply an external strain of $\ve_{22} = 0.06$ with a strain rate of $\dot{\ve}_{22} = 0.12$.

In both cases, as the number of cores increases from 64 to 2048, the total time to solution significantly decreases, demonstrating effective parallelization where the significant improvement from $64$ to $216$ is mainly due to the reduced sizes of the subdomains. However, the parallel efficiency shows a decline at higher core counts from $2\,048$ to $4\,096$ due to the huge size of the coarse space that reflects into communication overhead and load imbalance. Notably, the recycling option improves performance, reducing the preconditioner computation time by a significant factor, that contributes to a faster solution. While both configurations exhibit strong scalability, we notice that the recycling option does not affect the number of the iterations and does not increase substantially the solving time. This yields better parallel efficiency and reduced computational cost, making it a favorable choice for large-scale simulations.

\subsection{Weak Scalability}

\begin{table}[]
    \centering
    \begin{tabular}{ccccccccc}
        \hline
         	\multicolumn{9}{c}{No recycling}\\
        $Cores$ & $Coarse$ & $it_{GMRES}$ & $it_N$ & $T_{PC}$ & $T_{Sol}$ & $T_{Ass}$ & $T_{Tot}$ & $T_{N avg}$ \\
        \hline
        \hline
        64   & 1\,053  & 21.2 & 6.6  & 225   & 52    & 29.1   & 309.7 &  9.4  \\
        216  & 4\,115  & 23.4 & 7.2  & 279.4 & 70.3  & 31.3 & 385.4 & 10.7 \\
        512  & 10\,465 & 25.6 & 7.8  & 320.2 & 95.5  & 41.4   & 459.6 & 11.8 \\
        1728 & 37\,730 & 30.9 & 8.4  & 494.4 & 195.9 & 42.6   & 734.2 & 17.5 \\ 
        4096 & 92\,610 & 34.7 & 8.8  & 567.8 & 299.7 & 47.2   & 918.4 & 20.9 \\
        \hline
    \end{tabular}
    \centering
    \begin{tabular}{rrrrrrrrr}
        \hline
      	\multicolumn{9}{c}{With recycling}\\
        $Cores$ & $Coarse$ & $it_{GMRES}$ & $it_N$ & $T_{PC}$ & $T_{Sol}$ & $T_{Ass}$ & $T_{Tot}$ & $T_{N avg}$ \\
        \hline
        \hline
        64   & 1\,053  & 21.5 & 6.6  & 100.6  & 53.4  & 29.3 & 185.4 &  5.6  \\
        216  & 4\,115  & 23.9 & 7.2  & 114.1  & 71.9  & 30.7 & 219.1 &  6.1 \\
        512  & 10\,465 & 26.4 & 7.8  & 129.4  & 98.5  & 41.0 & 269.3 &  6.9 \\
        1728 & 37\,730 & 32.4 & 8.4  & 159.5  & 203.6 & 42.5 & 405.8 &  9.7 \\
        4096 & 92\,610 & 36.8 & 8.8  & 171.6  & 317.5 & 47.2 & 538.2 &  12.2 \\
        \hline
    \end{tabular}
    \caption{Weak scalability for GSDW*(T+R)-RGDSW without recycling option (top) and with recycling option (bottom). $16 \times 16 \times 16$ elements per subdomain. Computed the first 5 time iterations. All times are given in seconds.}
    \label{tab:weakRe}
\end{table}
With the same simulation set up of the \Cref{sec:strongCube}, in \Cref{tab:weakRe} we report a weak scalabitlity test again performed with GSDW*(T+R)-RGDSW coarse space. We keep fixed the subdomain sizes with $16 \times 16 \times 16$ elements and we analyze the performance when increasing proportionally with the number of cores. Both configurations exhibit a gradual increase in the computational time, with a more pronounced growth in the case without the recycling option.

Notably, the recycling option improves performance across all core counts, significantly reducing the preconditioner time. We note in fact, that from 64 to 4096 cores this increases by less than a factor of 2. Time for assembling the linear systems shows almost perfect weak scalability, while the time to compute the solution does not scale due to the increasing number of both, iteration per linear system and Newton iteration per time step.

For the reasons just described above, the average time per Newton iteration $T_{Navg}$ increase with the number of the subdomains. Despite this, the recycling strategy appears significantly faster, indicating improved efficiency in the iterative solution process. Considering the time per Newton iteration, we have a parallel efficiency of $46\%$ on $4\,096$ cores. This is an effect of slightly increasing \ac{GMRES} iteration counts per linear solve and an increase of the size of the coarse space. Considering instead the time per \ac{GMRES} iteration, we obtain an excellent parallel efficiency of $79\%$, which is satisfying due to the increasing coarse space size.

\begin{figure}[]
    \includegraphics[width=1\textwidth,trim={0 0 0 0},clip]{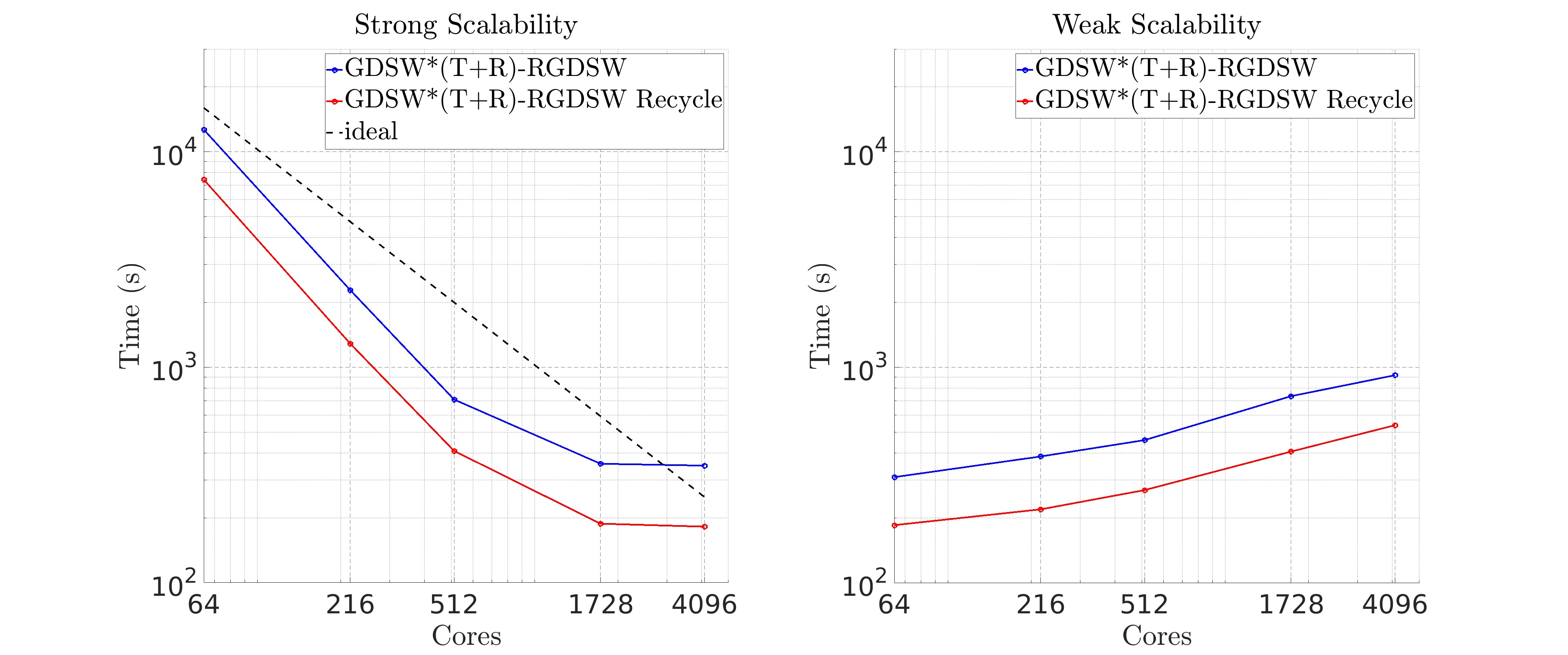}
{
  \caption{Strong and weak scalability for GDSW*-RGDSW with and without recycling option over 5 time iterations.}\label{fig:strongcube}
}
\end{figure}

\begin{figure}[]
    \includegraphics[width=1\textwidth,trim={0 150 0 200},clip]{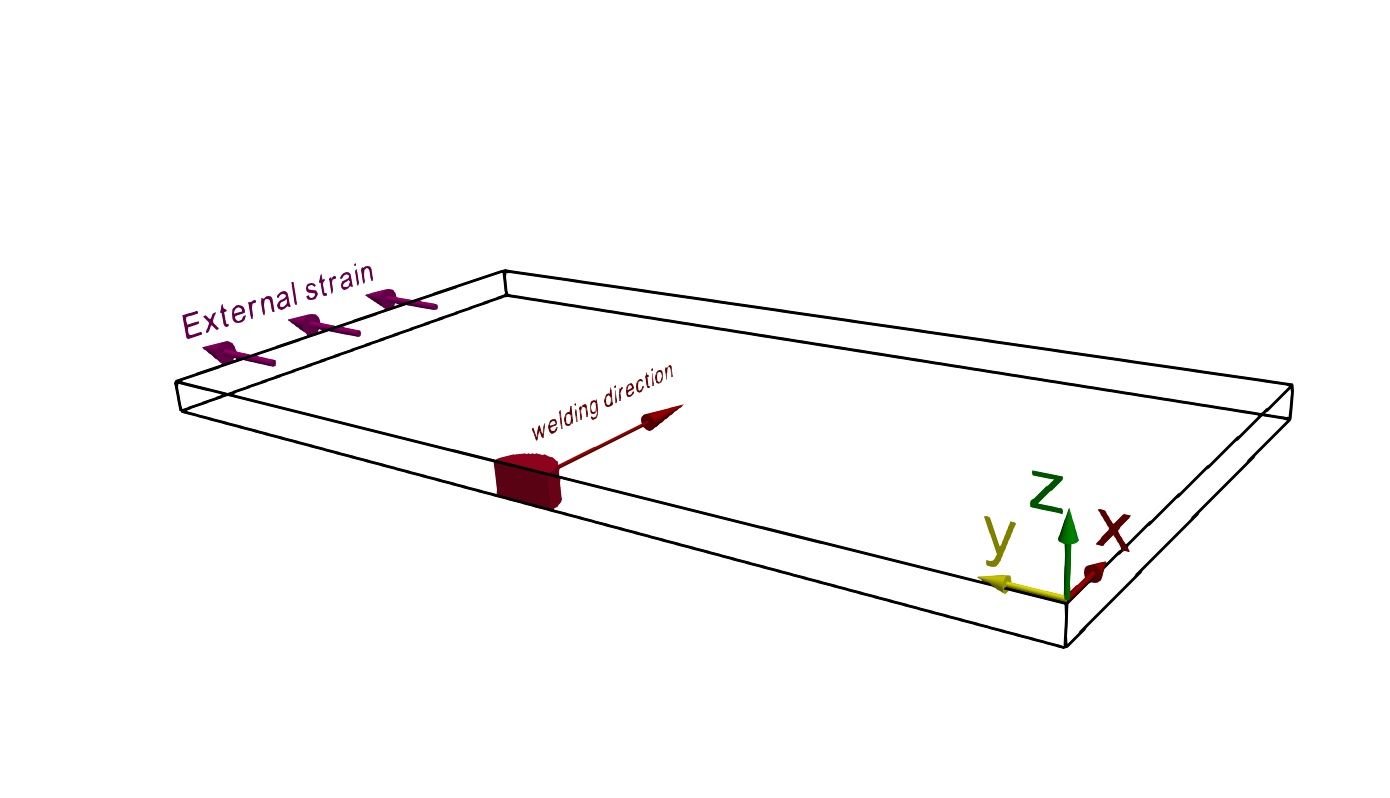}
  	\caption{Sketch of the thin plate geometry and displacements boundary conditions. In red the melting pool, initially placed with the center on the edge of the plate.}\label{fig:thintest}
\end{figure}

\subsection{2D domain decomposition}\label{sec:2Ddom}

\begin{table}[]
    \centering
    \begin{tabular}{ccccccc}
                \hline
        $PC\,Type$ & $Coarse$ & $it_{GMRES}$ & $it_{N}$ & $T_{PC}$ & $T_{Sol}$ & $T_{Tot}$  \\
         \hline
        \hline
          GDSW(T+R)-GDSW       & 19\,445 &  28.8  & 6.5  & 159.4 & 569.2    & 805.7 \\
          GDSW(T+R)-RGDSW      & 17\,478  & 28.2  & 6.5  & 147.5 & 501.2    & 722  \\
          RGDSW(T+R)-RGDSW     & 3\,777   & 294.3 & 6.5 & 114.7 & 3\,590   & 3\,779  \\ 
          \hline
          \hline
          GDSW(T)-GDSW       & 1\,1642  & 35.4  & 6.5 & 128.5 & 496.3  & 700.2 \\
          GDSW(T)-RGDSW      & 9\,681   & 35.4  & 6.5 & 120.3 & 459.6  & 654.7  \\
          RGDSW(T)-RGDSW     & 6\,615   & 302.9 & 6.5 & 107.5 & 3\,222 & 3\,405  \\ 
         \hline
    \end{tabular}
    \caption{ Solver performance with respect to different choices for the coarse space. $16\times64\times1$ subdomains with $8\times4\times9$ finite elements. The truncation tolerance is set to $1e-4$. All times are given in seconds.}
    \label{tab:comparison2D}
\end{table}
 To replicate the practical application of a \ac{LBW} process, we perform a series of numerical simulations on a thin plate using a 2D domain decomposition approach. The plate has dimensions $30 \times 15 \times 1$, and the simulation runs for $0.13$s real time. Initially,  during the first $0.1$s of the initialization phase the laser is positioned at the edge of the plate (see \Cref{fig:thintest}) and no external strain is applied. After this, the laser moves along the x-direction at a velocity of $16.67\,mm/s$ for $0.2$s, while an external strain is enforced with $\ve_{22} = 0.03$ and a strain rate of $\dot{\ve}_{22} =0.06$.

\Cref{tab:comparison2D} presents a comparison of different coarse spaces in this experiment, with the recycling option enabled. The results indicate that the GDSW-RGDSW approach, both with and without rotation, provides the best balance between computational cost and iteration count. As expected, the RGDSW(T)-RGDSW and RGDSW(T+R)-RGDSW configurations exhibit a significant increase in \ac{GMRES} iterations, leading to longer solving times. Overall,  even if GDSW(T)-RGDSW has slightly more \ac{GMRES} iterations than its rotational counterpart, it stands out as the most efficient configuration. Its superior performance, combined with a smaller coarse space that is half the size, makes it the optimal choice.\\
For this choice of the coarse space we perform two strong scalability tests reported in \Cref{tab:strongApp}. The two test are both computed with the same setup described just before using a finer discretization in space, namely $9\,492\,552$ and $30\,790\,760$ \ac{DOFs}. The smaller problem shows a perfect scalability up to $2\,048$ cores and then deteriorates when switching to $4\,096$ cores with a parallel efficiency of $55\%$, due to the coarse space that affect mainly the solution time. As showed in the same table, when moving to a higher number of \ac{DOFs}, we observe an excellent parallel efficiency of $76\%$ even with the larger number of subdomains. This is a well known phenomena in \ac{DDMs}, since strong scalability always deteriorates when the coarse space size gets large in comparison to the subdomain sizes. In fact, a direct solver application of the coarse level becomes prohibitive, due to the superlinear memory complexity of the sparse direct solvers. In such cases, on the coarse space can be recursively apllied a GDSW-type preconditioner, resulting in a three level method~\cite{heinlein2019three,heinlein2023multilevel}.

\begin{table}[]
    \centering
    \begin{tabular}{ccccc}
            \hline
         \multicolumn{3}{c}{9\,492\,552} \ac{DOFs} \\
        $Cores$ & $T_{tot} (T_{PC} + T_{Sol})$ & $it_{GMRES}$ & $T_{Tot}/\#Time it$ & $P_E$\\
        \hline
        \hline
        256  & 14\,846(6\,295+6\,740) & 49.7 & 114.2  & 1.00   \\
        512  & 6\,061(2\,371+2\,863)  & 40.9 & 46.6 & 1.20  \\
        1024 & 3\,099(939+1\,744)   & 44.7 & 23.8 & 1.19  \\
        2048 & 1\,812(470+1\,104)   & 39.2 & 13.9 & 1.01  \\
        4096 & 1\,700(344+1\,219)   & 42.2 & 13.1 & 0.55  \\
        \hline
        \hline
    \end{tabular}
    \hspace{1cm}
        \begin{tabular}{ccccc}
         \multicolumn{3}{c}{$30\,790\,760$} \ac{DOFs} \\
        $Cores$ & $T_{tot} (T_{PC} + T_{Sol})$ & $it_{GMRES}$ & $T_{Tot}/\#Time it$ & $P_E$\\
        \hline
        \hline
        512  & 10\,1647(22\,691+73\,668) & 188.5 & 781.9 & 1.00 \\
        1024 & 42\,789(8\,633+31\,642) & 137.0   & 329.1 & 1.18 \\
        2048 & 18\,739(3\,404+13\,985) & 114.7   & 144.1 & 1.35 \\
        4096 & 1\,390(1\,613+11\,532) & 113.2   & 107.0 & 0.91 \\
        \hline
    \end{tabular}
   \caption{Strong scalability for GDSW(T)-RGDSW with recycling option over 130 time iterations.}
    \label{tab:strongApp}
\end{table}

\subsection{Application to LBW simulation}\label{sec:Appl}
\begin{figure}[]
    \includegraphics[width=1\textwidth,trim={0 150 0 0},clip]{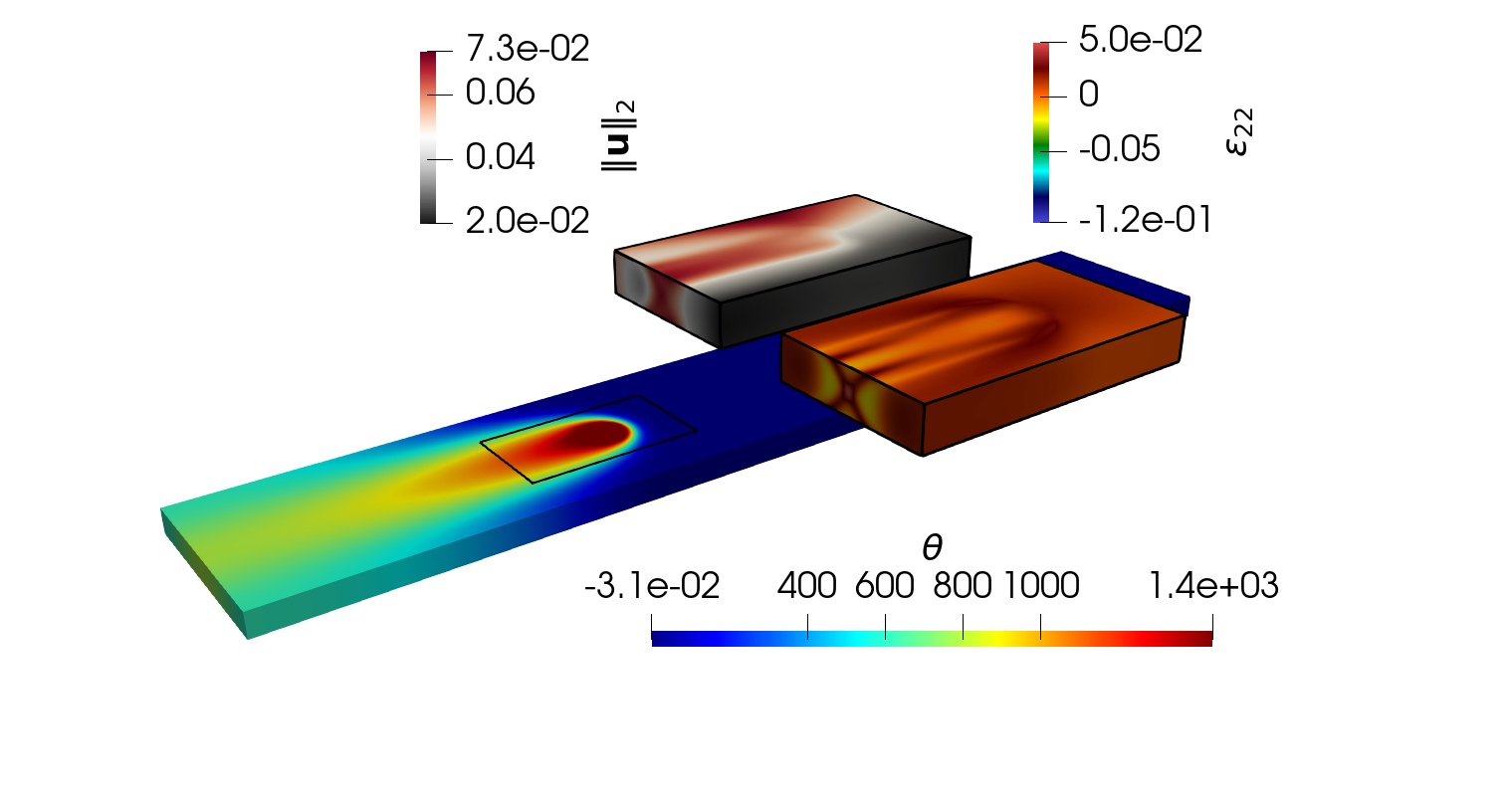}
  	\caption{Distribution of the temperature field $\theta$, the norm of the displacements $\Vert \bm{u} \Vert_2$ and the strain component into the direction of the external strain applied $\ve_{22}$ at the end of the simulation.}\label{fig:platesol}
\end{figure}

In this section we mimic a phYsically realistic simulation of a \ac{LBW} process. In particular we reproduce a simplified \ac{CTW} test to assess the weldability of thin plates in the \ac{LBW} process. We choose a plate of size $44.6 \times 8 \times 1 mm$ (see \Cref{fig:platesol}) and we split the domain into $128 \times 32 \times 1$ subdomain, each made of $7 \times 5 \times 20$ finite elements for a total of $12\,131\,028$ \ac{DOFs}. As before the laser is placed with the center on the edge of the plate, and no external strain is applied neither during the $0.1$s of the initialization phase nor during the following $0.7$s. The external strain is then applied for $0.3$s with $\ve_{22} = 0.03$ and $\dot{\ve}_{22} = 0.03$, leading to a total of $1.0$s real time until the end of the simulation. With a time step $\Delta t = 0.001$ this corresponds to $1\,000$ iterations. This simulation, performed with the solver previously described, i.e. \ac{GMRES} methods preconditioned with a two-level additive Schwarz preconditioner with GDSW(T)-RGDSW coarse space, with the recycling option active, and a truncation tolerance of $1e-4$. The coarse solver is performed on $64$ MPI ranks using cluster PARDISO. The simulation ran a total of $40\,684$s of runtime. This approximatively corresponds to $11$ hours and $20$ minutes of computation with $40.6$s per time iteration.

\section{Conclusions}
\label{sec:conclusions}

We have introduced different variants of monolithic GDSW coarse spaces for thermo-elastoplastic \ac{LBW} problems as well as a new and efficient parallel implementation in PETSc. Overall, the suggested preconditioners show a very high robustness and a good parallel scalability, which enables large scale simulations of fully-coupled thermo-mechanical simulations in the future. As a result of our numerical tests, the best candidate is a combination of GDSW or GDSW* for the mechanical problem with an RGDSW variant of reduced size for the temperature field.

\section*{Acknowledgments}
This project has received funding from the Deutsche For\-schungsgemeinschaft (DFG) as part of the Forschungsgruppe (Research Unit) 5134 ``Solidification Cracks During Laser Beam Welding -- High Performance Computing for High Performance Processing'' under DFG project number 434946896.
The authors gratefully acknowledge the scientific support and HPC resources provided by the Erlangen National High Performance Computing Center (NHR@FAU) of the Friedrich-Alexander-Universit\"at Erlangen-N\"urnberg (FAU) under the NHR project k109be10. NHR funding is provided by federal and Bavarian state authorities. NHR@FAU hardware is partially funded by the German Research Foundation (DFG) - 440719683. We would like to thank our project partners from Forschungsgruppe 5134 L.\ Scheunemann, J. Schr\"oder, and P. Hartwig for support with the interface to FEAP and providing the thermo-elastoplasticity formulation and all parameters in FEAP.

\bibliographystyle{siamplain}
\bibliography{literature}
\end{document}

%% file: ex_shared.tex
\usepackage{url}
\usepackage{color}
\usepackage{bm}
\usepackage{ulem}
\usepackage{comment}
\usepackage{siunitx}

\usepackage{siunitx}
\usepackage{lipsum}
\usepackage{acronym}
\usepackage{amsfonts}
\usepackage{graphicx}
\usepackage{epstopdf}
\usepackage{algorithmic}
\usepackage{multicol}
\usepackage{multirow} 
\ifpdf
  \DeclareGraphicsExtensions{.eps,.pdf,.png,.jpg}
\else
  \DeclareGraphicsExtensions{.eps}
\fi

\newcommand{\diff}{\,\mathrm{d}}
\newcommand{\ve         }{{\varepsilon }}
\newcommand{\Bve        }{{\bm{\varepsilon }}}
\newcommand{\Bsigma     }{{\bm{\sigma      }}}
\acrodef{CTW}{Controlled Tensile Weldability}
\acrodef{DOFs}{degrees of freedom}
\acrodef{DDMs}{domain decomposition methods}
\acrodef{GMRES}{Generalized Minimal Residual}
\acrodef{GDSW}{Generalized Dryja-Smith-Widlund}
\acrodef{LBW}{Laser Beam Welding}
\acrodef{PDEs}{partial differential equations}
\acrodef{PDE}{partial differential equation}

\newsiamremark{remark}{Remark}
\newsiamremark{hypothesis}{Hypothesis}
\crefname{hypothesis}{Hypothesis}{Hypotheses}
\newsiamthm{claim}{Claim}

\headers{Two-level Schwarz Preconditioners for laser beam welding}{T. Bevilacqua, A. Klawonn and M. Lanser}

\title{ Highly Scalable Two-level Monolithic Overlapping Schwarz Preconditioners for Thermo-elastoplastic laser beam welding problems\thanks{Submitted to the editors 03/2025.
\funding{This project has received funding from the Deutsche Forschungsgemeinschaft (DFG) as part of the Forschungsgruppe (Research Unit) 5134 ``Solidification Cracks During Laser Beam Welding -- High Performance Computing for High Performance Processing'' under DFG project number 434946896.}}}

\author{Tommaso Bevilacqua\thanks{Department of Mathematics and Computer Science, University of Cologne, Germany 
  (\email{tommaso.bevilacqua@uni-koeln.de, axel.klawonn@uni-koeln.de, martin.lanser@uni-koeln.de}).}
\and Axel Klawonn\footnotemark[2] \thanks{Center for Data and Simulation Science, University of Cologne, Germany}
\and Martin Lanser\footnotemark[2] \footnotemark[3]}

\usepackage{amsopn}
